\RequirePackage{ifpdf}
\ifpdf 
\documentclass[pdftex]{sigma}
\else
\documentclass{sigma}
\fi

\newcommand{\Om}{\Omega}
\newcommand{\pe}{\tau(e^{\perp})\leq \ep}
\newcommand{\ep}{\epsilon}
\newcommand{\om}{\omega}
\newcommand{\al}{\alpha}
\newcommand{\nm}{\Arrowvert}
\newcommand{\av}{\arrowvert}
\newcommand{\ii}{\infty}
\newcommand{\ol}{\overline}
\newcommand{\dt}{\delta}

\newcommand{\cl}{\mathcal{L}}
\newcommand{\im}{\Rightarrow}
\newcommand{\la}{\longrightarrow}

\newcommand{\gm}{\gamma}
\newcommand{\lb}{\lambda}

\begin{document}
\allowdisplaybreaks

\renewcommand{\PaperNumber}{023}

\FirstPageHeading

\ShortArticleName{A Banach Principle for Semif\/inite von Neumann
Algebras}

\ArticleName{A Banach Principle\\ for Semif\/inite von Neumann
Algebras}

\Author{Vladimir CHILIN~$^\dag$ and Semyon LITVINOV~$^\ddag$}

\AuthorNameForHeading{V. Chilin and S. Litvinov}

\Address{$^\dag$~Department of Mathematics, National University of
Uzbekistan, Tashkent 700095, Uzbekistan}
\EmailD{\href{mailto:chilin@ucd.uz}{chilin@ucd.uz}}
\Address{$^\ddag$~Department of Mathematics, Pennsylvania State University, \\
$\phantom{^\ddag}$~76 University Drive, Hazleton, PA 18202, USA}
\EmailD{\href{mailto:snl2@psu.edu}{snl2@psu.edu}}

\ArticleDates{Received November 25, 2005, in f\/inal form February
10, 2006; Published online February 20, 2006}

\Abstract{Utilizing the notion of uniform equicontinuity for
sequences of functions with the values in the space of measurable
operators, we present a non-commutative version of the Banach
Principle for $L^\infty$.} 


\Classification{46L51}

\section{Introduction}

Let $(\Om,\Sigma,\mu)$ be a probability space. Denote by $\cl
=\cl(\Om,\mu)$ the set of all (classes of) complex-valued
measurable functions on $\Om$. Let $\tau_{\mu}$ stand for the
measure topology in $\cl$. The classical Banach Principle may be
stated as follows.

\medskip

\noindent {\bf Classical Banach Principle.} {\it Let $(X, \nm
\cdot \nm)$ be a Banach space, and let $a_n:(X, \nm \cdot \nm)\to
(\cl,\tau_{\mu})$ be a sequence of continuous linear maps.
Consider the following properties:
\begin{enumerate}\vspace{-1mm}\itemsep=0pt

\item[\rm (I)] the sequence $\{ a_n(x)\}$ converges almost
everywhere (a.e.) for every $x\in X$;

\item[\rm (II)] $ a^{\star}(x)(\om)=\sup_n \av a_n(x)(\om)\av
<\ii$ a.e. for every  $x\in X$;

\item[\rm (III)] {\rm (II)} holds, and the maximal operator
$a^{\star}:(X, \nm \cdot \nm )\to (\cl, \tau_{\mu})$ is continuous
at $0$;

\item[\rm (IV)] the set $\{ x\in X: \{ a_n(x)\} \text {\ converges
a.e.} \}$ is closed in $X$.\vspace{-1mm}
\end{enumerate}

\noindent Implications {\rm (I)} $\im$ {\rm (II)} $\im$ {\rm
(III)} $\im$ {\rm (IV)} always hold. If, in addition, there is a
set $D\subset X$, $\ol D=X$, such that the sequence $\{ a_n(x) \}$
converges a.e.\ for every $x\in D$, then all four conditions {\rm
(I)--(IV)} are equivalent.}

\medskip

The Banach Principle is most often and successfully applied in the
context $X=(L^p, \nm \cdot \nm_p)$, $1\leq p<\ii$. At the same
moment, in the case $p=\ii$ the uniform topology in $L^{\ii}$
appears to be too strong for the ``classical'' Banach Principle to
be ef\/fective in $L^{\ii}$. For example, continuous functions are
not uniformly dense in $L^{\ii}$.

In \cite{BJ}, employing the fact that the unit ball $L_1^{\ii}=\{
x\in L^{\ii}: \nm x \nm_{\ii}\leq 1 \}$ is complete in
$\tau_{\mu}$, the authors suggest to consider the measure topology
in $L^{\ii}$ replacing $(X, \nm \cdot \nm )$ by
$(L_1^{\ii},\tau_{\mu})$. Note that, since $L_1^{\ii}$ is not a
linear space, geometrical complications occur, which in \cite{BJ}
are treated with the help of the following lemma.

\begin{lemma}\label{lemma 0.1}
{\it If $N(x,\dt)=\{ y\in L^{\ii}_1: \nm y-x\nm_1\leq \dt \}$,
$x\in L^{\ii}_1$, $\dt>0$, then $N(0,\dt)\subset
N(x,\dt)-N(x,\dt)$ for any $x\in L_1^{\ii}$, $\dt>0$.}
\end{lemma}

An application of the Baire category theorem yields the following replacement of ${\rm (I)} \im {\rm (II)}$.

\begin{theorem}[\cite{BJ}]\label{Theorem 0.2}   Let $a_n:L^{\ii}\to \cl$ be a sequence of
$\tau_{\mu}$-continuous linear maps such that the sequence $\{
a_n(x) \}$ converges a.e.\ for all $x\in L^{\ii}$. Then the
maximal operator $a^{\star}(x)(\om)=\sup_n \av a_n(x)(\om)\av$,
$x\in L^{\ii}$, is $\tau_{\mu}$-continuous at $0$ on $L_1^{\ii}$.
\end{theorem}

At the same time, as it is known~\cite{BJ}, even for a sequence $a_n:L^{\ii}\to
L^{\ii}$ of contractions, in which case condition~(II) is clearly
satisf\/ied, the maximal operator $a^{\star}:L_1^{\ii}\to
L_1^{\ii}$ may  be not $\tau_{\mu}$-continuous at $0$, i.e., 
(II) does not necessarily imply (III), whereas a replacement of the implication (III) $\im$
(IV) does hold:

\begin{theorem}[\cite{BJ}]\label{Theorem 0.3}  Assume that each $a_n:L^{\ii}\to \cl$ is linear, condition
{\rm (II)} holds with $X=L^{\ii}$, and the maximal operator
$a^{\star}:L^{\ii}\to \cl$ is $\tau_{\mu}$-continuous at $0$ on
$L_1^{\ii}$. Then the set $\{ x\in L_1^{\ii}: \{ a_n(x) \} \text
{\ converges a.e.} \}$ is closed in $(L_1^{\ii}, \tau_{\mu})$.
\end{theorem}

A non-commutative Banach Principle for measurable operators
af\/f\/iliated with a semif\/inite von Neumann algebra was
established in \cite{GL}. Then it was ref\/ined and applied in
\cite{LM,CLS,CL}. In \cite{CL} the notion of uniform
equicontinuity of a sequence of functions into $L(M,\tau)$ was
introduced. The aim of this study is to present a non-commutative
extension of the Banach Principle for $L^{\ii}$ that was suggested
in~\cite{BJ}. We were unable to prove a verbatim operator version
of Lemma~\ref{lemma 0.1}. Instead, we deal with the mentioned
geometrical obstacles via essentially non-commutative techniques,
which helps us to get rid of some restrictions in~\cite{BJ}.
First, proof of Lemma~\ref{lemma 0.1} essentially depends on the
assumption that the functions in $\cl$ be real-valued while the
argument of the present article does not employ this condition.
Also,
 our approach eliminates the assumption of the f\/initeness of measure.

\section{Preliminaries}

Let $M$ be a semif\/inite von Neumann algebra acting on a Hilbert
space $H$, and let $B(H)$ denote the algebra of all bounded linear
operators on $H$. A densely-def\/ined closed operator $x$ in $H$
is said to be {\it affiliated} with $M$ if $y^{\prime}x\subset
xy^{\prime}$ for every $y^{\prime}\in B(H)$ with $
y^{\prime}z=zy^{\prime}$, $z\in M $.
 We denote by $P(M)$ the complete lattice of all projections in $M$.
Let $\tau$ be a faithful normal semif\/inite trace on $M$. If $I$
is the identity of $M$, denote $e^{\perp}= I - e $, $e\in P(M)$.
An operator $x$ af\/f\/iliated with $M$ is said to be {\it
$\tau$-measurable} if for each $\ep >0$ there exists  a projection
$e\in P(M)$ with $\tau (e^{\perp})\leq \ep$ such that $eH$ lies in
the domain of the operator $x$.
 Let $L=L(M,\tau)$ stand for the set of all
$\tau$-measurable operators af\/f\/iliated with $M$. Denote $\nm
\cdot \nm$ the uniform norm in $B(H)$. If for any given $\ep>0$
and $\dt>0$ one sets
\[
 V(\ep, \dt)=\{x\in L: \nm xe\nm
\leq \dt \text { for some } e\in P(M) \text { with } \tau
(e^{\perp})\leq \ep \},
\]
then the topology $t_{\tau}$ in $L$ def\/ined by the family $\{
V(\ep, \dt): \ep>0, \dt>0 \}$ of neighborhoods of zero is called a
{\it measure topology}.

\begin{theorem}[\cite{Se}, see also \cite{Ne}]\label{Theorem 1.1}
 $(L,t_{\tau})$
is a complete metrizable topological $*$-algebra.
\end{theorem}

\begin{proposition}\label{Proposition 1.2} For any $d>0$, the sets $M_d=\{ x\in M: \nm x \nm \leq d \}$
and $M_d^h=\{ x\in M_d:x^*=x \}$ are $t_{\tau}$-complete.
\end{proposition}
\begin{proof} Because $(L,t_{\tau})$ is a complete metric space, it is enough to
show that $M_d$ and $M_d^h$ are (sequentially) closed in
$(L,t_{\tau})$.
 If $M_d\ni x_n\to_{t_{\tau}}x\in L$, then $0\leq x_n^*x_n\leq d\cdot I$ and,
due to Theorem \ref{Theorem 1.1}, $x_nx_n^*\to_{t_{\tau}}x^*x$.
Since $\{ x\in L:x\ge 0\}$ is $t_{\tau}$-complete, we have $0\leq
x^*x\leq d\cdot I$, which implies that $x\in M_d$. Therefore,
$M_d$ is closed in $(L,t_{\tau})$. Similarly, it can be checked
that $M_d^h$ is closed in $(L,t_{\tau})$.
\end{proof}

A sequence $\{ y_n \}\subset L$ is said to converge {\it almost
uniformly} (a.u.)  to $y\in L$ if for any given $\ep >0$ there
exists a projection $e\in P(M)$ with
 $\pe$ satisfying $\nm (y-y_n)e \nm \to 0$.

\begin{proposition}\label{Proposition 1.3}  If $\{ y_n\} \subset L$, then the conditions
\begin{enumerate}\vspace{-1mm}\itemsep=0pt
\item[\rm (i)] $\{y_n\}$ converges a.u. in $L$; \item[\rm (ii)]
for every $\ep >0$ there exists $e\in P(M)$ with $\pe$ such that
$\nm (y_m-y_n)e\nm \to 0$ as $m,n\to \ii$ \vspace{-1mm}
\end{enumerate}
are equivalent.
\end{proposition}

\begin{proof}
Implication ${\rm (i)}\im {\rm (ii)}$ is trivial. ${\rm (ii)}\im
{\rm (i)}$: Condition (ii) implies that the sequence $\{ y_n\}$ is
fundamental in measure. Therefore, by Theorem~\ref{Theorem 1.1},
one can f\/ind $y\in L$ such that $y_n\to y$ in~$t_{\tau}$. Fix
$\ep>0$, and let $p\in P(M)$ be such that $\tau(p^{\perp})\leq
\ep/2$ and $\nm (y_m-y_n)p \nm\to 0$ as $m,n\to \ii$. Because the
operators $ y_n$, $n\ge 1$, are measurable, it is possible to
construct such a projection $q\in P(M)$ with $\tau (q^{\perp})\leq
\ep/2$ that $\{ y_nq\} \subset M$. Def\/ining $e=p\land q$, we
obtain $\pe$, $y_ne=y_nqe\in M$, and
\[
\nm y_me-y_ne \nm = \nm (y_m-y_n)pe\nm \leq \nm (y_m-y_n)p\nm \la
0,
\]
$m,n\to \ii$. Thus, there exists $y(e)\in M$ satisfying $\nm
y_ne-y(e)\nm \to 0$. In particular, $y_ne\to y(e)$ in $t_{\tau}$.
On the other hand, $y_ne\to ye$ in $t_{\tau}$, which implies that
$y(e)=ye$. Hence, $\nm (y_n-y)e\nm \to 0$, i.e.\ $y_n\to y$ a.u.
\end{proof}

The following is a non-commutative Riesz's theorem \cite{Se}; see
also \cite{GL}.

\begin{theorem}\label{Theorem 1.4} {\it If $\{ y_n\} \subset L$ and $y=t_{\tau}-\lim\limits_{n\to \ii} y_n$,
then $y=a.u.-\lim\limits_{k\to \ii}y_{n_k}$ for some subsequence
$\{ y_{n_k} \} \subset \{ y_n \}$.}
\end{theorem}

\section[Uniform equicontinuity for sequences of maps into $L(M,\tau)$]{Uniform equicontinuity for sequences
of maps into $\boldsymbol{L(M,\tau)}$}

Let $E$ be any set. If $a_n:E\to L$, $x\in E$, and $b\in M$ are
such that $\{ a_n(x)b \} \subset M$, then we denote
\[
S(x,b)=S(\{ a_n \}, x, b)=\sup_n \nm a_n(x)b \nm.
\]
Def\/inition below is in part due to the following fact.

\begin{lemma}\label{Lemma 2.1}  Let $(X,+)$ be a semigroup, $a_n:X\to L$
be a sequence of additive maps. Assume that $\bar x\in X$ is such
that for every $\ep>0$ there exist a sequence $\{ x_k\} \subset X$
and a projection $p\in P(M)$ with $\tau (p^{\perp})\leq \ep$ such
that
\begin{enumerate}\vspace{-1mm}\itemsep=0pt
\item[\rm (i)] $\{ a_n(\bar x+x_k)\}$ converges a.u. as $n\to \ii$
for every $k$; \item[\rm (ii)] $S(x_k,p)\to 0$, $k\to \ii$.
\vspace{-1mm}
\end{enumerate}
Then the sequence $\{ a_n (\bar x)\}$ converges a.u.~in $L$.
\end{lemma}
\begin{proof}
Fix $\ep>0$, and let $\{ x_k\} \subset X$ and $p\in P(M)$, $\tau
(p^{\perp})\leq \ep/2$, be such that conditions~(i) and (ii) hold.
Pick $\dt>0$ and let $k_0=k_0(\dt)$ be such that $S(x_{k_0},p)\leq
\dt/3$. By Proposition~\ref{Proposition 1.3}, there is a
projection $q\in P(M)$ with $\tau(q^{\perp})\leq \ep/2$ and a
positive integer $N$ for which the inequality
\[
\nm (a_m(\bar x+x_{k_0})-a_n(\bar x+x_{k_0}))q\nm \leq \frac
{\dt}3
\]
holds whenever $m,n\ge N$. If one def\/ines $e=p\land q$, then
$\pe$ and
\begin{gather*}
\nm (a_m(\bar x)-a_n(\bar x))e \nm \leq \nm (a_m(\bar
x+x_{k_0})-a_n(\bar x+x_{k_0}))e\nm
\\
\qquad{}+ \nm a_m(x_{k_0})e\nm + \nm a_n(x_{k_0})e\nm \leq \dt
\end{gather*}
for all $m,n\ge N$. Therefore, by Proposition \ref{Proposition
1.3}, the sequence $\{ a_n(\bar x)\}$ converges a.u.\ in $L$.
\end{proof}

Let $(X,t)$ be a topological space, and let $a_n:X\to L$ and
$x_0\in X$ be such that $a_n(x_0)=0$, $n=1, 2, \dots$. Recall that
the sequence $\{ a_n \}$ is {\it equicontinuous} at $x_0$ if,
given $\ep>0$ and $ \dt>0$, there is a neighborhood $U$ of $x_0$
in $(X,t)$ such that $a_nU\subset V(\ep,\dt), n=1, 2, \dots$,
i.e., for every $x\in U$ and every $n$ one can f\/ind a projection
$e=e(x,n)\in P(M)$ with $\pe$ satisfying $\nm a_n(x)e \nm \leq
\dt$.

\medskip

\noindent {\bf Def\/inition.} Let $(X,t)$, $a_n:X\to L$, and
$x_0\in X$ be as above. Let $x_0\in E\subset X$. The sequence $\{
a_n \}$ will be called {\it uniformly equicontinuous} at $x_0$ on
$E$ if, given $\ep>0, \dt>0$, there is a neighborhood $U$ of $x_0$
in $(X,t)$ such that for every $x\in E\cap U$ there exists a
projection $e=e(x)\in P(M)$, $\pe$, satisfying $S(x,e) \leq \dt$.

\medskip

As it can be easily checked, the uniform equicontinuity is a
non-commutative generalization of the continuity of the maximal
operator, a number of equivalent forms of which are presented
in~\cite{BJ}.

 Let $\rho$ be an invariant metric in $L$
compatible with $t_{\tau}$ (see Theorem~\ref{Theorem 1.1}).

\begin{lemma}\label{Lemma 2.2} Let $d>0$. If a sequence $a_n:M\to L$ of additive maps is
uniformly equicontinuous at $0$ on $M_d^h$, then it is also
uniformly equicontinuous at $0$ on $M_d$.
\end{lemma}

\begin{proof}
Fix $\ep>0$, $\dt>0$. Let $\gm>0$ be such that, given $x\in
M_d^h$, $\rho(0,x)<\gm$, there is $e=e(x)\in P(M)$ for which $\tau
(e^{\perp})\leq \ep/2$ and $S(x,e)\leq \dt/2$ hold. Pick $x\in
M_d$ with $\rho(0,x)<\gm$. We have $x={\rm Re}\,(x)+i \,{\rm
Im}\,(x)$, where ${\rm Re}\,(x)= \frac {x+x^*}{2}$, ${\rm
Im}\,(x)=\frac {x-x^*}{2i}$. Clearly, ${\rm Re}\,(x),{\rm
Im}\,(x)\in M_d^h$ and $\rho(0,{\rm Re}\,(x))<\gm$, $\rho(0,{\rm
Im}\,(x))<\gm$. Therefore, one can f\/ind such $p,q\in P(M)$ with
$\tau(p^{\perp})\leq \ep/2$ and $\tau (q^{\perp})\leq \ep/2$ that
$S({\rm Re}\,(x),p)\leq \dt/2$ and $S({\rm Im}\,(x),q)\leq \dt/2$.
Def\/ining $r=p\land q$, we get $\tau(r^{\perp})\leq \ep$ and
\[
S(x,r)\leq S({\rm Re}\,(x),r)+S({\rm Im}\,(x),r)\leq S({\rm
Re}\,(x),p)+S({\rm Im}\,(x),q)\leq \dt,
\]
implying that the sequence $\{ a_n\}$ is uniformly equicontinuous
at $0$ on $M_d$.
\end{proof}

\begin{lemma}\label{Lemma 2.3} Let a sequence $a_n:M\to L$ of additive maps be uniformly
equicontinuous at $0$ on~$M_d$ for some $0<d\in \mathbb Q$. Then
$\{ a_n \}$ is also uniformly equicontinuous at $0$ on $M_s$ for
every $0<s\in \mathbb Q$.
\end{lemma}

\begin{proof} Pick $0<s\in \mathbb Q$, and let $r=d/s$. Given $\ep>0$, $\dt>0$,
one can present such $\gm>0$ that for every $x\in M_d$ with
$\rho(0,x)<\gm r$ there is a projection $e=e(x)\in P(M)$, $\pe$,
satisfying $S(x,e)\leq \dt r$. Since $a_n$ is additive and $d,s\in
\mathbb Q$, we have $a_n(rx)=ra_n(x)$. Also, $rx\in M_d$ and
$\rho(0,rx)<\gm r$ is equivalent to $x\in M_s$ and
$\rho(0,x)<\gm$. Thus, given $x\in M_s$ with $\rho(0,x)<\gm$, we
have
\[
\nm a_n(x)e \nm = \frac 1r \cdot \nm a_n(rx)e \nm \leq \dt,
\]
meaning that the sequence $\{ a_n \}$ is uniformly equicontinuous
at $0$ on $M_s$.
\end{proof}

\section{Main results}

Let $0\in E\subset M$. For a sequence of functions
$a_n:(M,t_{\tau})\to L$, consider the following conditions

\begin{enumerate}\itemsep=0pt\vspace{-1mm}
\item[] (CNV($E$)) almost uniform convergence of $\{ a_n(x) \}$
for every $x\in E$; \item[] (CNT($E$)) uniform equicontinuity at
$0$ on $E$; \item[] (CLS($E$)) closedness in $(E,t_{\tau})$ of the
set $C(E)=\{ x\in E: \{ a_n(x)\} \text {\ converges a.u.} \}$.
\vspace{-1mm}
\end{enumerate}
In this section we will study relationships among the conditions
(CNV($M_1$)), (CNT($M_1$)), and (CLS($M_1$)).

\medskip

\noindent {\bf Remarks.} 1. Following the classical scheme (see
Introduction), one more condition can be added to this list,
namely, a non-commutative counterpart of the existence of the
maximal operator, which can be stated as \cite{GL}:
\[
\mbox{(BND($E$)) given $x\in E$ and $\ep>0$, there is $e\in P(M)$,
$\pe$, with $S(x,e)<\ii$}.
\]
This condition can be called a {\it pointwise uniform boundedness}
of $\{ a_n\}$ on $E$. It can be easily verif\/ied that (CNV($E$))
implies (BND($E$)). But, as it was mentioned in Introduction, even
in the commutative setting, (BND($M_1$)) does not guarantee
(CNT($M_1$)).

2. If $a_n$ is additive for every $n$, then (CNV($M$)) follows
from (CNV($M_1$)).

3. If $E$ is closed in $(M,t_{\tau})$ (for instance, if $E=M_d$,
or $E=M_d^h$; see Proposition \ref{Proposition 1.2}), then
(CLS($E$)) is equivalent to the closedness of $C(E)$ in
$(L,t_{\tau})$.

\medskip

In order to show that (CNV($M_1$)) entails (CNT($M_1$)), we will
provide some auxiliary facts.

\begin{lemma}\label{Lemma 3.1}  For any $0\leq x\in L$ and $e\in P(M)$, $x\leq 2(exe+e^{\perp}xe^{\perp})$.
\end{lemma}
\begin{proof} If $a=e-e^{\perp}$, then $a^*=a$, which implies that
\[
0\leq axa=exe-exe^{\perp}-e^{\perp}xe+e^{\perp}xe^{\perp}.
\]
Therefore, $exe^{\perp}+e^{\perp}xe\leq exe+e^{\perp}xe^{\perp}$,
and we obtain
\begin{gather*}
x=(e+e^{\perp})x(e+e^{\perp})\leq
2(exe+e^{\perp}xe^{\perp}).\tag*{\qed}
\end{gather*} \renewcommand{\qed}{}
\end{proof}

For $y\in M$, denote $l(y)$ the projection on $\ol {yH}$, and let
$r(y)=I-n(y)$, where $n(y)$ denotes the projection on $\{ \xi \in
H: y\xi =0 \}$. It is easily checked that $l(y^*)=r(y)$, so, if
$y^*=y$, one can def\/ine $s(y)=l(y)=r(y)$. The projections
$l(y)$, $r(y)$, and $s(y)$ are called, respectively, a {\it left
support} of $y$, a {\it right support} of $y$, and a {\it support}
of $y=y^*$. It is well-known that $l(y)$ and $r(y)$ are equivalent
projections, in which case one writes $l(y)\sim r(y)$. In
particular, $\tau(l(y))=\tau(r(y))$, $y\in M$. If $y^*=y\in M$,
$y_+=\int_0^{\ii}\lb dE_{\lb}$, and $y_-=-\int_{-\ii}^0\lb
dE_{\lb}$, where $\{ E_{\lb}\}$ is the spectral family of $y$,
then we have $y=y_+-y_-$, $y_+=s(y_+)ys(y_+)$, and
$y_-=-s(y_+)^{\perp}ys(y_+)^{\perp}$.

\noindent The next lemma is, in a sense, a non-commutative
replacement of Lemma 0.1.

\begin{lemma}\label{Lemma 3.2} Let $y^*=y\in M$, $-I\leq y\leq I$. Denote $e_+=s(y_+)$.
If $x\in M$ is such that $0\leq x\leq I$, then
\[
-I\leq y-e_+xe_+\leq I \qquad \text {and}\qquad -I\leq
y+e_+^{\perp}xe_+^{\perp}\leq I.
\]
\end{lemma}
\begin{proof} Because $e_+xe_+\ge 0$, we have $y-e_+xe_+\leq y\leq I$; analogously,
$-I\leq y+e_+^{\perp}xe_+^{\perp}$. On the other hand, since we
obviously have $e_+xe_+\leq e_+$, $e_+^{\perp}xe_+^{\perp}\leq
e_+^{\perp}$, $e_+ye_+\leq e_+$, and $e_+^{\perp}ye_+^{\perp}\ge
-e_+^{\perp}$, one can write
\begin{gather*}
y-e_+xe_+=y_+-y_--e_+xe_+=y_++e_+^{\perp}ye_+^{\perp}-e_+xe_+ \ge
y_+-e_+^{\perp}-e_+=y_+-I\ge -I
\end{gather*}
and
\begin{gather*}
y+e_+^{\perp}xe_+^{\perp}=e_+ye_+-y_-+e_+^{\perp}xe_+^{\perp} \leq
e_+-y_-+e_+^{\perp}=I-y_-\leq I,
\end{gather*}
which f\/inishes the proof.
\end{proof}

\begin{lemma}\label{Lemma 3.3} $aV(\ep, \dt)b\subset V(2\ep,\dt)$ for
all $\ep>0$, $\dt>0$, and $a,b\in M_1$.
\end{lemma}

\begin{proof} Let $x\in V(\ep,\dt)$.  There exists $e\in P(M)$ such that
$\pe$ and $\nm xe\nm \leq \dt$. If we denote $q=n(e^{\perp}b)$,
then
\[
bq=(e+e^{\perp})bq=ebq+e^{\perp}bn(e^{\perp}b)=ebq.
\]
Besides, we have $q^{\perp}=r(e^{\perp}b)\sim l(e^{\perp}b)\leq
e^{\perp}$, which implies that $\tau(q^{\perp})\leq \ep$. Now, if
one def\/ines $p=e\land q$, then $\tau (p^{\perp})\leq 2\ep$ and
\[
\nm axbp \nm = \nm axbqp \nm =\nm axebqp\nm \leq \nm axeb \nm \leq
\nm a \nm \cdot \nm xe \nm \cdot \nm b \nm \leq \dt.
\]
Therefore, $axb\in V(2\ep,\dt)$.
\end{proof}
\begin{lemma}[\cite{GL}]\label{Lemma 3.4}   Let $f$ be the spectral projection of $b\in M$,
$0\leq b\leq I$, corresponding to the interval $[1/2, 1]$. Then
\begin{enumerate}\vspace{-1mm}\itemsep=0pt
\item[\rm (i)] $\tau(f^{\perp}) \leq 2\cdot \tau (I-b)$; \item[\rm
(ii)] $f=bc$ for some $c\in M$ with $0\leq c\leq 2\cdot I$.
\vspace{-1mm}
\end{enumerate}
\end{lemma}

We will also need the following fundamental result.

\begin{theorem}[\cite{Ka}] \label{Theorem 3.5}
 Let $a:M\to M$ be a positive linear map such that $a(I)\leq I$.
Then $a(x)^2\leq a(x^2)$ for every $x^*=x\in M$.
\end{theorem}

The next theorem represents a non-commutative extension of Theorem
\ref{Theorem 0.2}.

\begin{theorem}\label{Theorem 3.6}
Let $a_n:M\to L$ be a (CNV($M_1$)) sequence of positive
$t_{\tau}$-continuous linear maps such that $a_n(I)\leq I$,
$n=1,2, \dots $. Then the sequence $\{ a_n\}$ is also
(CNT($M_1$)).
\end{theorem}

\begin{proof}
Fix $\ep>0$ and $\dt>0$. For $N\in \mathbb N$ def\/ine
\[
F_N=
 \left \{
x\in M_1^h:
 \sup_{n\ge N} \nm (a_N(x)-a_n(x))b\nm \leq \dt
 \text {\ for some\ } b\in M, 0\leq b \leq I,
\tau (I-b)\leq \ep \right \}.
\]
Show that the set $F_N$ is closed in $(M_1^h,\rho)$.
 Let $\{ y_m\} \subset F_N$ and $\rho(y_m,\bar x)\to 0$ for some $\bar x \in L$.
 It follows from Proposition`\ref{Proposition 1.2} that $\bar x\in M_1^h$. We have
$a_1(y_m)\to a_1(\bar x)$ in $t_{\tau}$, which, by
Theorems~\ref{Theorem 1.1} and~\ref{Theorem 1.4},
 implies that there is
a subsequence $\{ y_m^{(1)} \} \subset \{ y_m \}$ such that
$a_1(y_m^{(1)})^*\to a_1(\bar x)^*$ a.u. Similarly, there is a
subsequence $\{ y_m^{(2)}\}\subset \{ y_m^{(1)} \}$ for which
$a_2(y_m^{(2)})^*\to a_2(\bar x)^*$ a.u. Repeating this process
and def\/ining $x_m=y_m^{(m)}\in F_N$, $m=1,2, \dots$, we obtain
\[
a_n(x_m)^*\la a_n(\bar x)^* \ \ \text {a.u.,} \quad  m\to \infty ,
\quad n=1, 2, \dots .
\]
By def\/inition of $F_N$, there exists a sequence
$\{b_{m}\}\subset M$, $0\leq b_{m}\leq I$, $\tau (I-b_{m})\leq
\ep$, such that $\sup_{n\ge N} \nm (a_N(x_m)-a_n(x_m))b_m \nm \leq
\dt$ for every $m$.
 Because $M_1$ is weakly compact,
there are a subnet $\{b_{\al}\} \subset \{b_{m}\}$ and  $b\in M$
 such that $b_{\al}\to b$ weakly, i.e. $(b_{\al}\xi,\xi)\to
(b\xi,\xi)$ for all $\xi \in H$. Clearly $0\leq b\leq I$. Besides,
by the well-known inequality (see, for example \cite{BR}),
\[
\tau (I-b)\leq \liminf_{\al} \tau (I-b_{\al}) \leq \ep.
\]
 We shall show that $\sup_{n\ge N} \nm (a_N(\bar x)-a_n(\bar x))b\nm \leq \dt$.
Fix $n\ge N$. Since $ a_{k}(x_{m})^*\to a_{k}(\bar x)^*$ a.u.,
$k=n,N$, given $\sigma >0$, there exists a projection $e\in P(M)$
with $\tau (e^{\perp})\leq \sigma$ satisfying
\[
\left \| e\left (a_k(x_m)-a_k(\bar x)\right ) \right \|= \left \|
\left (a_k(x_m)^* -a_k(\bar x)^*\right )e\right \| \la 0, \qquad
m\to \ii, \quad k=n,N.
\]
 Show f\/irst that
$\nm e(a_N(\bar x)-a_n(\bar x))b\nm \leq \dt$. For every $\xi,\eta
\in H$ we have
\begin{gather}
\left \av \left (e((a_N(x_m)-a_n(x_m))b_m-(a_N(\bar x)-a_n(\bar
x))b)\xi,\eta \right) \right \av
\nonumber\\
\qquad{}\leq \left \av \left (e(a_N(x_m)-a_n(x_m)-a_N(\bar
x)+a_n(\bar x))b_m\xi,\eta \right )
\right \av \nonumber\\
\qquad\phantom{{}\leq{}}{} +\left \av \left ((b_m-b)\xi,(a_N(\bar
x)^* -a_n(\bar x)^*)e\eta \right ) \right \av.\label{eq1}
\end{gather}
Fix $\gm >0$ and choose $m_0$ be such that
\begin{gather}\label{eq2}
\left \| e\left (a_k(x_m)-a_k(\bar x)\right )\right \| <\gm,
\qquad k=n,N
\end{gather}
whenever $m\ge m_0$. Since $b_{\al}\to b$ weakly, one can f\/ind
such an index $\al(\gm)$ that
\begin{gather}
\left \av \left ((b_{\al}-b)\xi,(a_N(\bar x)^*-a_n(\bar x)^*)e\eta
\right ) \right \av <\gm\label{eq3}
\end{gather}
 as soon as $\al \ge \al(\gm)$. Because $\{ b_{\al}\}$ is a subnet of $\{ b_m\}$,
 there is such an index $\al(m_0)$ that
$\{ b_{\al} \}_{\al \ge \al(m_0)} \subset \{ b_m\}_{m\ge m_0}$. In
particular, if $\al_0 \ge \max \ \{ \al(\gm), \al(m_0)\}$, then
$b_{\al_0}=b_{m_1}$ for some $m_1 \ge m_0$. It follows now from
\eqref{eq1}--\eqref{eq3} that, for all $\xi, \eta \in H$ with $\nm
\xi \nm =\nm \eta \nm =1$, we have
\begin{gather*}
\left \av \left (e(a_N(\bar x)-a_n(\bar x))b\xi,\eta \right )
\right \av
\leq \left \av \left (e(a_N(x_{m_1})-a_n(x_{m_1}))b_{m_1}\xi,\eta \right ) \right \av \\
\qquad{} +\left \av \left (e(a_N(x_{m_1})-a_n(x_{m_1})-a_N(\bar
x)+a_n(\bar x))b_{m_1}\xi,\eta
 \right ) \right \av \\
\qquad{} + \left \av \left ((b_{m_1}-b)
\xi,(a_N(\bar x)^*-a_n(\bar x)^*)e\eta \right ) \right \av \\
\qquad{} \leq \dt + \left \| e\left (a_N(x_{m_1}) -a_N(\bar
x)\right )\right \| +\nm e(a_n(x_{m_1})-a_n(\bar x)) \nm +\gm <
\dt +3\gm.
\end{gather*}
Due to the arbitrariness of $\gm>0$, we get
\[
\left \| e(a_N(\bar x)-a_n(\bar x))b\right \|=\sup_{\| \xi \| = \|
\eta \| =1}\left \av \left (e(a_N(\bar x)-a_n(\bar x))b\xi, \eta
\right ) \right \av \leq \dt.
\]
Next, we choose $e_{j}\in P(M)$ such that $\tau (e_{j}^{\perp})
\leq \frac {1}{j}$ and
\[
\left \| e_j\left (a_k(x_m)-a_k(\bar x)\right )\right \| \la 0
\text { \ as  \ } m\to \infty,\qquad k=n,N; \quad j=1,2, \dots .
\]
Since $e_j\to I$ weakly, $e_j(a_N(\bar x)-a_n(\bar x))b \to
(a_N(\bar x)-a_n(\bar x))b$ weakly, therefore,
\[
 \nm (a_N(\bar x)-a_n(\bar x))b\nm
 \leq \limsup_{j\to \ii}
\nm e_j(a_N(\bar x)-a_n(\bar x))b \nm \leq \dt.
\]
Thus, for every $n\ge N$ the inequality $\nm (a_N(\bar x)-a_n(\bar
x))b\nm \leq \dt$ holds, which implies that $\bar x\in F_N$ and
$\ol {F_N}=F_N$.

Further, as $\{ a_n(x)\}$ converges a.u.\ for every $x\in M_1$,
taking into account Proposition \ref{Proposition 1.3}, we obtain
\[
M_1^h=\bigcup^{\infty}_{N=1}F_N.
\]

 By Proposition \ref{Proposition 1.2}, the metric space $(M_1^h,\rho)$ is complete.
Therefore, using the Baire category theorem, one can present such
$N_0$ that $F_{N_0}$ contains an open set. In other words, there
exist $x_0\in F_{N_0}$ and $\gm_0\ge 0$ such that for any $x\in
M_1^h$ with $\rho(x_0,x)<\gm_0$ it is possible to f\/ind $b_x\in
M$, $0\leq b_x\leq I$, satisfying $\tau(I-b_x)\leq \ep$ and
\[
\sup_{n\ge N_0}\nm (a_{N_0}(x)-a_n(x))b_x \nm \leq \dt.
\]
Let $f_x$ be the spectral projection of $b_x$ corresponding to the
interval $[1/2, 1]$. Then, according to Lemma~\ref{Lemma 3.4},
$\tau(f_x^{\perp})\leq 2\ep$ and
\[
\sup_{n\ge N_0} \nm (a_{n_0}(x)-a_n(x))f_x\nm \leq 2\dt
\]
 whenever $x\in M_1^h$ and $\rho(x_0,x)<\gm_0$.
Since the multiplication in $L$ is continuous with respect to the
measure topology, Lemma~\ref{Lemma 3.3} allows us to choose
$0<\gm_1<\gm_0$ in such a way that $\rho(0,x)<\gm_1$ would imply
$\rho(0,ax^2b)<\gm_0$ for every $a,b\in M_1$. Denote
$e_+=s(x_0^+)$. Because $a_i:(M,\rho)\to (L,t_{\tau})$ is
continuous for each $i$, there exists such $0<\gm_2<\gm_1$ that,
given $x\in M$ with $\rho(0,x)<\gm_2$, it is possible to f\/ind
such a projection $p\in P(M)$, $\tau(p^{\perp})\leq \ep$, that
\[
\nm a_i(e_+x^2e_+)p\nm \leq \dt\qquad \text {and} \qquad \nm
a_i(e_+^{\perp}x^2e_+^{\perp})p\nm \leq \dt,
\]
$i=1, \dots , N_0$. Let $x\in M_1^h$ be such that
$\rho(0,x)<\gm_2$. Since $0\leq x^2\leq I$, Lemma~\ref{Lemma 3.2}
yields
\[
-I\leq x_0-e_+x^2e_+\leq I \qquad \text {and}\qquad  -I\leq
x_0+e_+^{\perp}x^2e_+^{\perp}\leq I,
\]
so, we have
\[
y=x_0-e_+x^2e_+\in M_1^h \qquad \text{and} \qquad
z=x_0+e_+^{\perp}x^2e_+^{\perp}\in M_1^h.
\]
Besides, $\rho(x_0,y)=\rho(0,-e_+x^2e_+)<\gm_0$, which implies
that there is $f_1\in P(M)$ such that $\tau(f_1^{\perp})\leq 2\ep$
and
\[
\sup_{n\ge N_0} \nm (a_{N_0}(y)-a_n(y))f_1\nm \leq 2\dt.
\]
Analogously, one f\/inds $f_2\in P(M)$, $\tau(f_2^{\perp})\leq
2\ep$, satisfying
\[
\sup_{n\ge N_0} \nm (a_{N_0}(z)-a_n(z))f_2\nm \leq 2\dt.
\]
As $\rho(0,x)<\gm_2$, there is $p\in P(M)$ with
$\tau(p^{\perp})\leq \ep$ such that the inequalities
\[
\nm a_i(e_+x^2e_+)p\nm \leq \dt\qquad \text{and}\qquad \nm
a_i(e_+^{\perp}x^2e_+^{\perp})p\nm \leq \dt
\]
hold for all $i=1, \dots , N_0$. Let $e=f_{x_0}\land f_1\land
f_2\land p$. Then we have $\tau(e^{\perp})\leq 7\ep$ and, for
$n>N_0$,
\begin{gather*}
\nm a_n(e_+x^2e_+)e\nm \leq \nm (a_{N_0}(x_0-e_+x^2e_+)-a_n(x_0-e_+x^2e_+) \\
\qquad{} +a_n(x_0)-a_{N_0}(x_0)+a_{N_0}(e_+x^2e_+))e\nm \leq \nm (a_{N_0}(y)-a_n(y))f_1e\nm \\
\qquad{} +\nm (a_{N_0}(x_0)-a_n(x_0))f_{x_0}e\nm +\nm
a_{N_0}(e_+x^2e_+)pe\nm \leq 5\dt.
\end{gather*}
At the same time, if $n\in \{ 1, \dots , N_0 \}$, then $\nm
a_n(e_+x^2e_+)e\nm = \nm a_n(e_+x^2e_+)pe\nm \leq \dt$, so
\[
\nm a_n(e_+x^2e_+)e\nm \leq 5\dt, \qquad n=1,2, \dots .
\]
Analogously,
\[
\nm a_n(e_+^{\perp}x^2e_+^{\perp})e\nm \leq 5\dt, \qquad n=1,2,
\dots .
\]
Next, by Lemma \ref{Lemma 3.1}, we can write $0\leq x^2\leq
2(e_+x^2e_++e_+^{\perp}x^2e_+^{\perp})$. Since $a_n$ is positive
for every~$n$, applying Theorem \ref{Theorem 3.5}, we obtain
\[
0\leq ea_n(x)^2e\leq ea_n(x^2)e\leq
2(ea_n(e_+x^2e_+)e+ea_n(e_+^{\perp}x^2e_+^{\perp})e).
\]
Therefore,
\[
\nm a_n(x)e \nm^2=\nm ea_n(x)^2e \nm \leq 20\dt, \qquad n=1,2,
\dots .
\]
Summarizing, given $\ep>0$, $\dt>0$, it is possible to f\/ind such
$\gm>0$ that for every $x\in M_1^h$ with $\rho(0,x)<\gm$ there is
a projection $e=e(x)\in P(M)$ such that $\tau(\ep^{\perp})\leq
7\ep$ and
\[
S(x,e)=\sup_n\nm a_n(x)e\nm \leq \sqrt {20\dt}.
\]
Thus, the sequence $\{a_n \}$ is (CNT($M_1^h$)), hence, by
Lemma~\ref{Lemma 2.2}, (CNT($M_1$)).
\end{proof}

Now we shall present a non-commutative extension of Theorem
\ref{Theorem 0.3}.

\begin{theorem}\label{Theorem 3.7} A {\rm (CNT($M_1$))} sequence $a_n:M\to L$ of additive maps is also
{\rm (CLS($M_1$))}.
\end{theorem}

\begin{proof} Let $\bar x$ belong to the $t_{\tau}$-closure of $C(M_1)$.
By Proposition~\ref{Proposition 1.2}, $\bar x\in M_1$. Fix
$\ep>0$. Since, by Lemma~\ref{Lemma 2.3}, the sequence $\{ a_n\}$
is (CNT($M_2$)), for every $k\in \mathbb N$, there is $\gm_k>0$
such that, given $x\in M_2$ with $\rho(0,x)<\gm_k$, one can f\/ind
a projection $p_k=p_k(x)\in P(M)$, $\tau(p_k^{\perp})\leq
\ep/2^k$, satisfying $S(x,p_k)\leq 1/k$. Let a sequence $\{ y_n\}
\subset C(M_1)$ be such that $\rho(\bar x, y_k)< \gm_k$. If we set
$x_k=y_k-\bar x$, then $x_k\in M_2$, $\rho (0,x_k)=\rho(\bar
x,x_k+\bar x)= \rho(\bar x,y_k)<\gm_k$, and $\bar x+x_k=y_k\in
C(M_1)$, $k=1,2, \dots$. If $e_k=p_k(x_k)$, then
$\tau(e_k^{\perp})\leq \ep/2^k$ and also $S(x_k,e_k)\leq 1/k$.
Def\/ining $e=\land_{k=1}^{\ii}$, we obtain $\pe$ and
$S(x_k,e)\leq 1/k$. Therefore, by Lemma~\ref{Lemma 2.1}, the
sequence $\{ a_n(\bar x)\}$ converges a.u., i.e. $\bar x\in
C(M_1)$.
\end{proof}

The following is an immediate consequence of the previous results
of this section.

\begin{theorem}\label{Theorem 3.8} Let $a_n:M\to L$ be a sequence of positive $t_{\tau}$-continuous
linear maps such that $a_n(I)\leq I$, $n=1,2, \dots$. If $\{
a_n\}$ is $({\rm CNV}(D))$ with $D$ being $t_{\tau}$-dense in $M_1$,
then conditions $({\rm CNV}(M_1))$, $({\rm CNT}(M_1))$, and $({\rm CLS}(M_1))$ are
equivalent.
\end{theorem}

\section{Conclusion}

First we would like to stress that, due to Theorem \ref{Theorem 3.6}, 
when establishing the almost uniform convergence of a
sequence $\{a_n(x)\}$ for all $x\in L^{\ii}(M,\tau) = M$, the
uniform equicontinuity at~0 on $M_1$ of the sequence $\{a_n\}$
is assumed. Also, as it is noticed in \cite{BJ}, the above
formulation is important because,
 for example, if $\{a_n\}$ are bounded operators in a
 non-commutative $L^p$-space, $1\leq p < \ii$, one may want to show that not
only do these operators fail to converge a.u., but they fail so
badly that $\{a_n\}$ may fail to converge a.u.\ on any class of
operators which is $t_\tau$-dense in $M$.

\subsection*{Acknowledgements}
S. Litvinov is partially supported by the 2004 PSU RD Grant.

\LastPageEnding

\end{document}